\documentclass[11pt]{article}
\usepackage{amssymb}
\usepackage{latexsym}
\usepackage{amsfonts, amsmath}
\usepackage{epsfig}
\include{epsfig.sty}

\newcommand{\demo}{ {\it   Proof. }}
\newcommand{\qed}{$\Box$\bigskip}

\newtheorem{thm}{Theorem}[section]
\newtheorem{prop}[thm]{Proposition}
\newtheorem{lem}[thm]{Lemma}
\newtheorem{cor}[thm]{Corollary}

\newtheorem{defn}[thm]{Definition}
\newtheorem{rmk}[thm]{Remark}

\begin{document}

\title{The automorphism group of a free-by-cyclic group in rank 2}

\author{O.\ Bogopolski \\ \small{Institute of Mathematics of SBRAS,}\\ {\small Novosibirsk, Russia} \\ \small{and Centre
de Recerca Matem\`{a}tica} \\ \small{e-mail: groups@math.nsc.ru} \\ \\
\\ A.\ Martino \\ \small{Dept.\ Mat.\ Apl.\ IV, Univ.\ Pol.\ Catalunya, (Barcelona, Spain)}
\\ \small{and Centre de Recerca Matem\`{a}tica} \\ \small{e-mail: Armando.Martino@upc.edu}  \\ \\ \\ E.\ Ventura \\ \small{Dept.\ Mat.\ Apl.\
III, Univ.\ Pol.\ Catalunya, (Barcelona, Spain)}
\\ \small{and Centre de Recerca Matem\`{a}tica} \\ \small{e-mail: enric.ventura@upc.edu} }

\maketitle

\begin{abstract}
Let $\phi$ be an automorphism of a free group $F_n$ of rank $n$, and
let $M_{\phi}=F_n \rtimes_{\phi} \mathbb{Z}$ be the corresponding
mapping torus of $\phi$. We study the group $Out(M_{\phi})$ under
certain technical conditions on $\phi$. Moreover, in the case of
rank 2, we classify the cases when this group is finite or virtually
cyclic, depending on the conjugacy class of the image of $\phi$ in
$GL_2(\mathbb{Z})$.
\end{abstract}

\section{Introduction}

Let $F_n$ be the free group of rank $n$ freely generated by
$x_1,\ldots ,x_n$, and let us denote automorphisms $\phi \in
Aut(F_n)$ as acting on the right, $x\mapsto x\phi$. In this paper we
consider extensions of finitely generated free groups by the
infinite cyclic group (\emph{[f.g.~free]-by-$\mathbb{Z}$} groups,
for short). More precisely, for any given $\phi \in Aut(F_n)$, we
consider the mapping torus, $M_{\phi}=F\rtimes_{\phi} \mathbb{Z}$,
of $\phi$ i.e. the extension of $F_n$ presented by
 $$
M_{\phi}=\langle x_1,\dots ,x_n ,\, t \mid t^{-1}x_it=x_i\phi \,\,\,
(i=1,\ldots ,n)\rangle.
 $$
The aim of the paper is to study the outer automorphism group of
such groups, $Out(M_{\phi})$. We shall give partial results for
arbitrary rank $n$, and a complete description for the cases
$n=1,2$.

To help avoiding possible confusions, we will use Greek letters
(such as $\phi$ or $\psi$) to denote automorphisms of $F_n$, and
capital Greek letters (such as $\Phi$ or $\Psi$) to denote
automorphisms of $M_{\phi}$. Accordingly, for every word $w\in F_n$,
we shall write $\gamma_w$ to denote the inner automorphism of $F_n$
given by right conjugation by $w$, $x\gamma_w =w^{-1}xw$. And, for
every element $g\in M_{\phi}$, we shall write $\Gamma_g$ to denote
the inner automorphism of $M_{\phi}$ given by right conjugation by
$g$, $x\Gamma_g =g^{-1}xg$. As usual, $Inn(G)$ denotes the group of
inner automorphisms of a group $G$, and $Out(G)=Aut(G)/Inn(G)$.

Although [f.g.~free]-by-$\mathbb{Z}$ groups have received a great
deal of attention in recent years, there has been no real systematic
study of their automorphisms. In fact, it still seems to be an open
question whether or not they have finitely generated or finitely
presented automorphism groups. Having said that, there are certain
cases in which the automorphism group is understood. For instance,
when $M_{\phi}$ is word hyperbolic, it is known to have finite outer
automorphism group (this can be deduced from papers \cite{BF,Br}).
However, note that in the rank 2 case, the group $M_{\phi}$ can
never be hyperbolic. In fact, by a result of Nielsen
(see~\cite[proposition~5.1]{LS}), $([x_1,x_2])\phi$ and so
$t^{-1}[x_1,x_2]t$, must be conjugate to $[x_1,x_2]^{\pm 1}$ in
$F_2$. Hence, $M_{\phi}$ contains a free abelian subgroup of rank 2
implying that $M_{\phi}$ is not hyperbolic.

The case of rank 1 is straightforward to analyse. If $\phi \in
Aut(\mathbb{Z})$ is the identity then $M_{\phi}=F_1 \rtimes_{\phi}
\mathbb{Z} =\mathbb{Z}\times \mathbb{Z}$ and so,
$Out(M_{\phi})=Aut(M_{\phi})=GL_2(\mathbb{Z})$. Otherwise, $\phi$ is
the inversion and it is easy to check that $Out(M_{\phi})$ is finite
in this case.

%
%

The paper is organised as follows. In Section~\ref{n}, we analyse
$Aut(M_{\phi})$ for arbitrary $n$ but under certain technical
restrictions for $\phi$ (see Theorems~\ref{evector} and~\ref{prod}
below). After dedicating Section~\ref{back2} to recall a standard
classification of $2\times 2$ matrices, the main result comes in
Section~\ref{2}, where we analyse $Aut(M_{\phi})$ in the case when
the underlying free group has rank $n=2$, and without conditions on
$\phi$. The rank 2 case is doubtless the easiest to deal with, but
we believe that some of our methods may be of general interest. For
instance, our detailed look at the parabolic case is certain to be
of use in the more general UPG case (the definition of UPG
automorphisms can be found in \cite{BFH}). The information obtained
in Section~\ref{2} is summarised in the subsequent theorem, which is
the main result of the paper.

\begin{thm}\label{main}
Let $F_2 =\langle a,b\rangle$ be a free group of rank 2, let $\phi
\in Aut(F_2)$, and consider the mapping torus $M_{\phi}=F_2
\rtimes_{\phi} \mathbb{Z}$. Let $\phi^{\rm \, ab} \in GL_2
(\mathbb{Z})$ be the map induced by $\phi$ on $F_2^{\rm \, ab}\cong
\mathbb{Z}^2$ (written in row form with respect to $\{a,b\}$).
\begin{itemize}
\item[i)] If $\phi^{\rm \, ab}=I_2$, then $Out(M_{\phi})\cong (\mathbb{Z}^2 \rtimes C_2)\rtimes GL_2(\mathbb{Z})$, where
$C_2$ is the cyclic group of order 2 acting on $\mathbb{Z}^2$ by
sending $u$ to $-u$, $u\in \mathbb{Z}^2$, and where
$GL_2(\mathbb{Z})$ acts trivially on $C_2$ and naturally on
$\mathbb{Z}^2$ (thinking vectors as columns there).
\item[ii)] If $\phi^{\rm \, ab}=-I_2$, then $Out(M_{\phi}) \cong PGL_2(\mathbb{Z}) \times C_2$.
\item[iii)] If $\phi^{\rm \, ab}\neq -I_2$ and does not have 1 as an eigenvalue then $Out(M_{\phi})$ is finite.
\item[iv)] If $\phi^{\rm \, ab}$ is conjugate to $\left( \smallmatrix 1 & k \\ 0 & -1 \endsmallmatrix \right)$ for
some integer $k$, then $Out(M_{\phi})$ has an infinite cyclic
subgroup of finite index.
\item[v)] If $\phi^{\rm \, ab}$ is conjugate to $\left(\smallmatrix 1 & k \\ 0 & 1 \\ \endsmallmatrix
\right)$ for some integer $k\neq 0$, then $Out(M_{\phi})$ has an
infinite cyclic subgroup of finite index.
\end{itemize}
Furthermore, for every $\phi \in Aut(F_2)$, $\phi^{\rm \, ab}$ fits
into exactly one of the above cases.
\end{thm}

We note that in rank 2, $M_{\phi}$ is always the fundamental group of a 3-manifold, namely the corresponding mapping torus of a once punctured torus,
say $M$. When $\phi^{\rm \, ab}$ is hyperbolic, $M$ is a hyperbolic manifold and it is known that, in this case, by Mostow rigidity, $Out(\pi_1
M)=Out(M_{\phi})$ is finite. This fact is contained in Theorem~\ref{main}~(iii), which we prove by elementary methods. Although this is the generic
case, we believe there is value in analysing the entire situation, especially since the proof of Theorem~\ref{main}~(v) is rather more involved than
one might expect.


We remark that the isomorphism type of $M_{\phi}$ depends only on
the conjugacy class of the outer automorphism $[\phi]^{\pm 1}\in
Out(F_n) =Aut(F_n)/Inn(F_n)$ determined by $\phi^{\pm 1}$ (this is
straightforward to verify, see~Lemma~\ref{semi} below). As an
application of our main result, we prove that this characterises the
isomorphism class of $M_{\phi}$, when $n=2$.

\begin{cor}\label{rk2}
Let $F_2 =\langle a,b\rangle$ be a free group of rank 2 and let
$\phi,\, \psi \in Aut(F_2)$. The groups $M_{\phi}$ and $M_{\psi}$
are isomorphic if and only if $[\phi]$ and $[\psi]^{\pm 1}$ are
conjugate in $Out(F_2)$.
\end{cor}

\demo The if part is straightforward and proven in Lemma~\ref{semi}
below.

Let us assume $M_{\phi}\cong M_{\psi}$ (and so, $Out(M_{\phi})\cong
Out(M_{\psi})$). By Theorem~\ref{main}, $\phi$ and $\psi$ fit each
one into exactly one of the items (i) to (v). But the groups
$(\mathbb{Z}^2 \rtimes C_2)\rtimes GL_2(\mathbb{Z})$,
$PGL_2(\mathbb{Z}) \times C_2$, a finite group, and a
virtually-$\mathbb{Z}$ group are not pairwise isomorphic (note, for
example, that the first one contains $\mathbb{Z}^2$ and the second
one is virtually free). So, $\phi$ and $\psi$ fit simultaneously
into case (i), or case (ii), or case (iii), or cases (iv)-(v). In
order to distinguish also between cases (iv) and (v), note that
$\left( \smallmatrix 1 & k \\ 0 & -1
\endsmallmatrix \right)^2 =\left( \smallmatrix 1 & 0 \\ 0 & 1 \endsmallmatrix \right)$ while $\left(\smallmatrix 1 & k
\\ 0 & 1 \\ \endsmallmatrix \right)$ has infinite order; Lemma~\ref{centre} then, says that $M_{\phi}\cong M_{\psi}$ has
non-trivial and trivial centre, respectively.

Suppose $\phi$ and $\psi$ both fit into case (i) or both into case
(ii). Then $\phi^{\rm \, ab}=\psi^{\rm \, ab}$, and so $\phi
\psi^{-1}$ is an inner automorphism of $F_2$. Thus, $[\phi ]=[\psi
]$ in $Out(F_2)$.

Suppose $\phi$ and $\psi$ both fit into case (iii) and so,
$\phi^{\rm \, ab}$ and $\psi^{\rm \, ab}$ do not have 1 as an
eigenvalue. By Theorem~\ref{evector}, $F_2$ is the unique normal
subgroup of $M_{\phi}$ and of $M_{\psi}$ with quotient isomorphic to
$\mathbb{Z}$. Hence, the given isomorphism $\Xi \colon M_{\phi} \to
M_{\psi}$ must preserve the free part, $F_2$, and so must send $t$
to $t^{\epsilon}w$, for some $\epsilon =\pm 1$ and $w\in F_2$. Let
$\xi \in Aut(F_2)$ be the restriction of $\Xi$ to $F_2$. For every
$x\in F_2$, the relation $t^{-1}xt =x\phi$ (valid in $M_{\phi}$)
implies $w^{-1}(x\xi\psi^{\epsilon})w=w^{-1}t^{-\epsilon}(x\xi
)t^{\epsilon}w =x\phi\xi$ in $M_{\psi}$. So, we have $\xi
\psi^{\epsilon}\gamma_w =\phi \xi$ in $Aut(F_2)$, which means that
$[\phi]$ and $[\psi]^{\epsilon}$ are conjugate in $Out(F_2)$.

Suppose $\phi$ and $\psi$ both fit into case (iv), say $\phi^{\rm \,
ab}$ is conjugate to $\left( \smallmatrix 1 & k_1
\\ 0 & -1 \endsmallmatrix \right)$ and $\psi^{\rm \, ab}$ is conjugate to $\left( \smallmatrix 1 & k_2 \\ 0 & -1
\endsmallmatrix \right)$. It is straightforward to verify that these two matrices are conjugate to each other in
$GL_2(\mathbb{Z})$ if and only if $k_1$ and $k_2$ have the same
parity. So, $\phi^{\rm \, ab}$ (and similarly $\psi^{\rm \, ab}$) is
conjugate to either $\left( \smallmatrix 1 & 0 \\ 0 & -1
\endsmallmatrix \right)$ or $\left( \smallmatrix 1 & 1 \\ 0 & -1
\endsmallmatrix \right)$. Accordingly, $M_{\phi}$ (and $M_{\psi}$)
is isomorphic to either
 $$
M_{_{\left( \smallmatrix 1 & 0 \\ 0 & -1 \endsmallmatrix \right)}}
=\langle a,b,t \mid t^{-1}at=a,\, t^{-1}bt=b^{-1}\rangle
 $$
or
 $$
M_{_{\left( \smallmatrix 1 & 1 \\ 0 & -1 \endsmallmatrix \right)}}
=\langle a,b,t \mid t^{-1}at=ab,\, t^{-1}bt=b^{-1}\rangle.
 $$
But, looking at their abelianisations, these two groups are not
isomorphic to each other (the first abelianises to $\mathbb{Z}^2
\oplus \mathbb{Z}/2\mathbb{Z}$ and the second to $\mathbb{Z}^2$).
Thus, $k_1$ and $k_2$ have the same parity and so, $[\phi]$ and
$[\psi]$ are conjugate in $Out(F_2)$.

Finally, suppose $\phi$ and $\psi$ both fit into case (v), say
$\phi^{\rm \, ab}$ is conjugate to $\left( \smallmatrix 1 & k_1 \\ 0
& 1
\endsmallmatrix \right)$ and $\psi^{\rm \, ab}$ is conjugate to $\left( \smallmatrix 1 & k_2 \\ 0 & 1
\endsmallmatrix \right)$. Easy computations show that $M_{\phi}^{\rm \, ab} \cong \mathbb{Z}^2 \oplus \mathbb{Z}/\vert
k_1 \vert \mathbb{Z}$ and $M_{\psi}^{\rm \, ab} \cong \mathbb{Z}^2
\oplus \mathbb{Z}/\vert k_2 \vert \mathbb{Z}$. Hence, $k_1 =\pm
k_2$. And, for every $k$, the matrix $\left( \smallmatrix 1 & 0 \\ 0
& -1 \endsmallmatrix \right)$ conjugates $\left( \smallmatrix 1 & k
\\ 0 & 1 \endsmallmatrix \right)$ into $\left( \smallmatrix 1 & -k
\\ 0 & 1
\endsmallmatrix \right)$. So, $[\phi]$ and $[\psi]$ are conjugate in $Out(F_2)$. \qed

\begin{rmk} \emph{
We remark that some results similar to Corollary~\ref{rk2} were
already known in the abelian context. For example, in the first
appendix of~\cite{GS} the authors show that two
$\mathbb{Z}^2$-by-cyclic groups are isomorphic if and only if the
corresponding actions (i.e. matrices in $GL_2(\mathbb{Z})$) are
conjugated one to the other or one to the inverse of the other (our
Corollary~\ref{rk2} states exactly the same changing $\mathbb{Z}^2$
to $F_2$). Also, W. Dicks extended that result to dimension 3 with
the extra hypothesis that the corresponding $3\times 3$ matrices
have no non-trivial eigenvectors of eigenvalue one (see section 4.5
of~\cite{H}).
 }\end{rmk}

For completeness, let us emphasise that Corollary~\ref{rk2} is no longer true for general rank. However, finding two explicit non-conjugate (neither
inverse-conjugate) outer automorphisms $\phi$ and $\psi$ of a free group such that $M_{\phi} \cong M_{\psi}$ does not seem to be straightforward. We
present the following example, which is the simplest in a series of examples we have found. We thank Joan Porti for pointing out to us certain
manifolds that allowed us to construct an explicit example in rank 6. We also thank Warren Dicks, who then provided a rank 3 example, and allowed us
to include it here. Also, we thank R. Vikent'ev (Bogopolski's student) for proposing some other similar examples in rank 3.

Dicks' example is the following. Consider the group presented by
 $$
G=\langle \, s,t \mid t^{-3}st^2st^{-1}s^{-1}ts^{-2}ts \, \rangle
 $$
On one hand we have
 $$
\begin{array}{rl}
G & \cong \langle \, s_0, s_1, s_2, s_3, t \mid s_1 =t^{-1}s_0 t,\,
s_2 =t^{-1}s_1 t,\, s_3 =t^{-1}s_2 t,\, s_3 s_1 s_2^{-1}s_1^{-2}s_0
=1 \, \rangle
\\ & \cong \langle \, s_0, s_1, s_2, t \mid t^{-1}s_0 t =s_1,\, t^{-1}s_1 t=s_2,\, t^{-1}s_2 t =s_0^{-1} s_1^2 s_2 s_1^{-1} \, \rangle
\\ & \cong M_{\phi},
\end{array}
 $$
where $\phi$ is the automorphism of $F_3 =\langle a,b,c\rangle$
given by $a\to b\to c\to a^{-1}b^2cb^{-1}$. On the other hand,
 $$
\begin{array}{rl}
G & \cong \langle \, t_{-2}, t_{-1}, t_0, t_1, s \mid t_{-1}
=s^{-1}t_{-2} s,\, t_0 =s^{-1}t_{-1} s,\, t_1 =s^{-1}t_0 s,\,
t_0^{-3} t_{-1}^2 t_{-2}^{-1}t_{-1}t_1 =1 \, \rangle \\ & \cong
\langle \, t_{-2}, t_{-1}, t_0, s \mid s^{-1}t_{-2} s =t_{-1},\,
s^{-1}t_{-1} s=t_0,\, s^{-1}t_0 s =t_{-1}^{-1} t_{-2} t_{-1}^{-2}
t_0^3 \, \rangle
\\ & \cong M_{\psi},
\end{array}
 $$
where $\psi$ is the automorphism of $F_3 =\langle a,b,c\rangle$
given by $a\to b\to c\to b^{-1}ab^{-2}c^3$.

This proves that $M_{\phi}\cong M_{\psi}$, while the outer
automorphisms $[\phi], \, [\psi] \in Out(F_3)$ are neither conjugate
nor inverse-conjugate since their abelianisations are two $3\times
3$ matrices with determinant -1 and 1, respectively. This concludes
the desired example.

\section{Results for general rank}\label{n}

Using the defining relations of $M_{\phi}$ under the form
$wt=t(w\phi)$ and $wt^{-1} =t^{-1}(w\phi^{-1})$, $w\in F_n$, it is
clear that in any element of $M_{\phi}$ we can push all the $t$
letters to one side. That is, we have a normal form in $M_{\phi}$
whereby we can write elements uniquely in the form $t^k w$, where
$k$ is an integer and $w\in F_n$.

Our first observation is that the isomorphism type of $M_{\phi}$
depends only on the outer automorphism determined by $\phi^{\pm 1}$,
up to conjugacy in $Out(F_n)$.

\begin{lem}\label{semi}
Let $F_n$ be a free group of rank $n$, let $\phi,\psi \in Aut(F_n)$
and consider
 $$
M_{\phi}=\langle x_1,\dots ,x_n ,\, t \mid t^{-1}x_it=x_i\phi \,\,\,
(i=1,\ldots ,n)\rangle
 $$
and
 $$
M_{\psi}=\langle x_1,\dots ,x_n ,\, s \mid s^{-1}x_is=x_i\psi \,\,\,
(i=1,\ldots ,n)\rangle.
 $$
If the automorphisms $\phi$ and $\psi$ are conjugate or
conjugate-inverse to each other in $Out(F_n)$, then $M_{\phi}$ and
$M_{\psi}$ are isomorphic. More precisely, if $\chi \in Aut(F_n)$ is
such that $\chi^{-1} \phi \chi = \psi^{\epsilon} \gamma_w$ for some
$\epsilon = \pm 1$ and $w\in F_n$, then the map $\Omega \colon
M_{\phi} \to M_{\psi}$, $x_i \mapsto x_i \chi$, $t\mapsto
s^{\epsilon}w$ extends to an isomorphism.
\end{lem}

\demo The map $\Omega$ is extended multiplicatively and one needs to
show that it is well defined. In order to do this, it is enough to
show that the relators in $M_{\phi}$ are all sent to the trivial
element in $M_{\psi}$. In fact, for every $i=1,\ldots ,n$, we have
 $$
\begin{array}{rcl}
t^{-1} x_i t (x_i \phi)^{-1} & \mapsto & w^{-1} s^{-\epsilon} (x_i
\chi) s^{\epsilon} w (x_i \phi \chi)^{-1} \\ & = & w^{-1} (x_i \chi
\psi^{\epsilon}) w (x_i \phi \chi)^{-1} \\ & = & (x_i \chi
\psi^{\epsilon} \gamma_w) (x_i \phi \chi)^{-1} \\ & = & (x_i \phi
\chi) (x_i \phi \chi)^{-1} \\ & = & 1\in M_{\psi}.
\end{array}
 $$
Thus, we have a well defined homomorphism $\Omega \colon M_{\phi}
\to M_{\psi}$, which is surjective by inspection. Moreover, if an
element $t^k x\in M_{\phi}$ is in the kernel of $\Omega$, we
immediately deduce that $k=0$ and $x\chi =1$. Since $\chi$ is an
automorphism of $F_n$, we have $x=1$ and $\Omega$ has trivial
kernel. $\Box$

We continue with some basic facts about the automorphism group of
$M_{\phi}$.

\begin{lem}\label{centre}
Let $F_n$ be a free group of rank $n$, and let $\phi \in Aut(F_n)$.
The group $M_{\phi}=F_n \rtimes_{\phi} \mathbb{Z}$ has non-trivial
centre if and only if $\phi^k=\gamma_w$ for some $k\neq 0$ and some
$w\in F_n$. If this equation holds and $n\geqslant 2$, then $w\phi
=w$.
\end{lem}

\demo The result is clear for $n=0,1$. So, we may assume that
$n\geqslant 2$.

A straightforward calculation shows that the element $t^kw^{-1} \in
M_{\phi}$ commutes with every $x\in F_n$ if and only if $\phi^k
=\gamma_w$. Similarly, $t^kw^{-1}$ commutes with $t$ if and only if
$w\phi =w$. So, $t^kw^{-1}$ is central in $M_{\phi}$ if and only if
$\phi^k =\gamma_w$ and $w\phi =w$. Hence, $M_{\phi}$ has non-trivial
centre if and only if $\phi^k =\gamma_w$ and $w\phi =w$ for some
integer $k$ and some $w\in F_n$ such that $(k,w)\neq (0,1)$.

Now, using the fact $n\geqslant 2$, we can simplify this. Assume the
equation $\phi^k =\gamma_w$ holds. Since, $\gamma_w =\phi^{-1}
\gamma_w \phi= \gamma_{w\phi}$ we have $w\phi =w$. Also, note that
$k=0$ implies $w=1$ (because $F_n$ has trivial centre for
$n\geqslant 2$). Thus, $M_{\phi}$ has non-trivial centre if and only
if $\phi^k =\gamma_w$ for some integer $k\neq 0$ and some $w\in
F_n$. \qed


Let $\Psi \in Aut(M_{\phi})$ and suppose it leaves $F_n$ invariant.
In this situation, its restriction to $F_n$, $\psi
=\Psi\vert_{F_n}$, is an endomorphism of $F_n$ that will be seen in
the next proposition to be always an automorphism. On the other
hand, factorising by the normal and $\Psi$-invariant subgroup $F_n$,
we get an automorphism $\overline{\Psi}$ of $M_{\phi}/F_n \cong
\mathbb{Z}$. If $\overline{\Psi}=Id$ we shall say that $\Psi$ is a
\emph{positive} automorphism of $M_{\phi}$. Otherwise,
$\overline{\Psi}$ is the inversion of $\mathbb{Z}$ and we say that
$\Psi$ is \emph{negative}. In any case, $t\Psi =t^{\epsilon}w$ for
some $w\in F$, where $\epsilon =\pm 1$ is the \emph{signum} of
$\Psi$.

\begin{prop}\label{inv}
Let $F_n$ be a free group of rank $n$, let $\phi \in Aut(F_n)$ and
consider $M_{\phi}=F_n \rtimes_{\phi} \mathbb{Z}$. Let $\Psi \in
Aut(M_{\phi})$ be such that $F_n \Psi \leqslant F_n$, and denote by
$\psi \colon F_n \to F_n$ its restriction to $F_n$. Then,
\begin{itemize}
\item[i)] $\psi$ is an automorphism of $F_n$,
\item[ii)] there exists $w\in F_n$ such that $\phi \psi = \psi \phi^{\epsilon} \gamma_w$, where $\epsilon$ is the signum of
$\Psi$. Furthermore, if $n\geqslant 2$ then the word $w$ is unique
and satisfies the equation $t\Psi=t^{\epsilon}w$.
\end{itemize}
\end{prop}

\demo Since $F_n\Psi \leqslant F_n$, we must have that $t\Psi =
t^{\pm 1} w$ for some $w\in F_n$ (otherwise, $t$ would not be in the
image of $\Psi$). Now, clearly, $F_n\Psi$ is a normal subgroup of
$M_{\phi}=M_{\phi}\Psi =\langle F_n\psi, t^{\pm 1}w \rangle$. Hence
any element $g\in M_{\phi}$ can be written in the form $g=(v\psi
)(t^{\pm 1}w)^k$ for some $v\in F_n$ and $k\in \mathbb{Z}$. And here
$g\in F_n$ if and only if $k=0$. Thus, $F_n \Psi=F_n$ and $\Psi$
induces an automorphism on $F_n$. This proves (i).

Let $\epsilon =\pm 1$ be the signum of $\Psi$, that is, $t\Psi =
t^{\epsilon}w$ for some $w\in F_n$. Applying $\Psi$ to both sides of
the equality $x\phi=t^{-1}xt$ we get
 $$
x\phi \psi =w^{-1}t^{-\epsilon}(x\psi )t^{\epsilon}w =x\psi
\phi^{\epsilon} \gamma_w,
 $$
for all $x\in F_n$. Hence, $\phi \psi = \psi \phi^{\epsilon}
\gamma_w$. Furthermore, if $n\geqslant 2$, this last equation can
only hold for a unique $w\in F_n$. This proves (ii).  $\Box$

In the following result, we impose certain hypothesis on $\phi$ to
ensure that every automorphism of $M_{\phi}$ leaves the free
subgroup invariant. Under these circumstances, computing
$Out(M_{\phi})$ is fairly straightforward.

\begin{thm}\label{evector}
Let $F_n$ be a free group of rank $n$, let $\phi \in Aut(F_n)$ and
consider $M_{\phi}=F_n \rtimes_{\phi} \mathbb{Z}$. Let
$M_{\phi}^{\rm \, ab}$ denote the abelianisation of $M_{\phi}$, and
$F_n^{\rm \, ab}\cong \mathbb{Z}^n$ the abelianisation of $F_n$
(which is not in general the image of $F_n \leqslant M_{\phi}$ in
$M_{\phi}^{\rm \, ab}$). Let $\phi^{\rm \, ab} \in GL_n
(\mathbb{Z})$ be the map induced by $\phi$ on $F_n^{\rm \, ab}$, and
$[\phi]$ be the class of $\phi$ in $Out(F_n)$. The following are
equivalent:
\begin{itemize}
\item[a)] $M_{\phi}^{\rm \, ab}$ is the direct sum of an infinite cyclic group and a finite abelian group,
\item[b)] the matrix $\phi^{\rm \, ab}$ does not have eigenvalue 1.
\item[c)] $F_n \leqslant M_{\phi}$ is the unique normal subgroup of $M_{\phi}$ with quotient isomorphic to $\mathbb{Z}$.
\end{itemize}

Furthermore, if these conditions hold then every automorphism of
$M_{\phi}$ leaves $F_n$ invariant,
 $$
Aut^+ (M_{\phi})=\{ \Psi \in Aut(M_{\phi}) \mid \Psi \mbox{ is
positive} \}
 $$
is a normal subgroup of $Aut(M_{\phi})$ of index at most 2 and
moreover, if $n\geqslant 2$, its image $Out^+(M_{\phi})$ in
$Out(M_{\phi})$ is also normal, of index at most two, and isomorphic
to $C([\phi])/\langle [\phi ]\rangle$, where $C([\phi])$ denotes the
centraliser of $[\phi]$ in $Out(F_n )$.
\end{thm}

\demo  To prove the equivalence of (a) and (b), note that
$M_{\phi}^{\rm \, ab}\cong \langle t \mid \,\, \rangle \oplus
F_n^{\rm \, ab}/ Im(\phi^{\rm \, ab}-Id)$. Then, $F_n^{\rm \,
ab}/Im(\phi^{\rm \, ab}-Id)$ is finite if and only if ${\mbox\rm{
rank}}_{_{\mathbb{Z}}}(Im(\phi^{\rm \, ab}-Id))={\mbox\rm{
rank}}_{_{\mathbb{Z}}}(F_n^{\rm \, ab})=n$. And this happens if and
only if ${\mbox\rm{ rank}}_{_{\mathbb{Z}}}(Ker(\phi^{\rm \,
ab}-Id))=0$, which is the same as saying that $\phi^{\rm \, ab}$ has
no eigenvalue 1. On the other hand, every epimorphism from
$M_{\phi}$ onto $\mathbb{Z}$ factors through $M_{\phi}^{\rm \, ab}$.
So, the equivalence between (a) and (c) is clear ($F_n$ being the
preimage in $M_{\phi}$ of the torsion subgroup of $M_{\phi}^{\rm \,
ab}$).

We shall now prove the remaining assertions under the assumption
that these three conditions hold. Consider the abelianisation map
$M_{\phi}\rightarrow M_{\phi}^{\rm \, ab}$. Since $F_n^{\rm \,
ab}/Im(\phi^{\rm \, ab}-Id)$ is the torsion subgroup of
$M_{\phi}^{\rm \, ab}$, its full pre-image $F_n$ is characteristic
in $M_{\phi}$. Hence, every automorphism of $M_{\phi}$ leaves $F_n$
invariant. At this point, it is clear that $Aut^+(M_{\phi})$ is a
normal subgroup of $Aut(M_{\phi})$ of index at most 2, and so is
$Out^+(M_{\phi})$ in $Out(M_{\phi})$.

Assuming $n\geqslant 2$, it remains to prove that
$Out^+(M_{\phi})\cong C([\phi])/\langle [\phi]\rangle$. Define the
map
 $$
\begin{array}{ccl} f \colon Aut^{+}(M_{\phi}) & \rightarrow & C([\phi])/\langle [\phi]\rangle \\ \Psi & \mapsto &
[\Psi|_{F_n}]\langle [\phi]\rangle. \end{array}
 $$
Note that by Proposition~\ref{inv}, the image of this map lies in
$C([\phi])/\langle [\phi]\rangle$. Clearly, $f$ is a homomorphism.

First, we will prove that $Im \, f=C([\phi])/\langle [\phi]
\rangle$. Let $\psi \in Aut(F_n)$ be such that $[\psi]\in
C([\phi])$. Then, $\phi\psi =\psi \phi \gamma_w$ for some $w\in
F_n$. In this situation, it is straightforward to verify that $\psi$
extends to a well defined automorphism $\Psi$ of $M_{\phi}$ by just
sending $t$ to $tw$. Clearly, $\Psi\in Aut^+(M_{\phi})$ and its
$f$-image is $[\psi]\langle [\phi]\rangle$.

Now we will prove that $Ker\, f =Inn(M_\phi)$. For every element
$g=t^k w\in M_{\phi}$, we have that $\Gamma_g
|_{F_n}=\phi^k\gamma_w$ and so, $\Gamma_g$ maps under $f$ to the
identity element of $C([\phi])/\langle [\phi]\rangle$. This means
that $Ker\, f \geqslant Inn(M_\phi)$. Conversely, let $\Psi\in
Ker\,f$. Then, $\Psi|_{F_n}=\phi^k\gamma_w =\Gamma_{t^kw} |_{F_n}$
for some integer $k$ and some $w\in F_n$. Also, applying $\Psi$ to
both sides of the equation $t^{-1}xt=x\phi$ and using the positivity
of $\Psi$ and the fact $n\geqslant 2$, it is straightforward to
check that $t\Psi =t(w\phi)^{-1}w =t\Gamma_{t^kw}$. Thus, $\Psi
=\Gamma_{t^k w}$. This completes the proof that $Ker\, f
=Inn(M_\phi)$ and so, $Out^+(M_{\phi})\cong C([\phi])/\langle
[\phi]\rangle$. $\Box$

The extreme opposite case to the one considered above is when $\phi$
is the identity automorphism (or, in fact, an inner automorphism).
We also calculate the automorphism group in this case.

\begin{thm}\label{prod}
Let $F_n =\langle x_1,\ldots ,x_n \rangle$ be a free group of rank
$n\geqslant 2$, and let $M=M_{Id}=F_n \times \mathbb{Z}$. Consider
the group $\mathbb{Z}^n \rtimes C_2$ where $C_2$ is the cyclic group
of order 2 which acts by sending $u$ to $-u$ for all $u\in
\mathbb{Z}^n$ (think $u$ as a column vector); also, consider the
action of $Aut(F_n)$ (and also $Out(F_n)$) on it given by the
trivial action on $C_2$, and the natural action after abelianisation
on $\mathbb{Z}^n$. Then, $Aut(M)\cong (\mathbb{Z}^n \rtimes
C_2)\rtimes Aut(F_n)$ and $Out(M)\cong (\mathbb{Z}^n \rtimes
C_2)\rtimes Out(F_n)$.
\end{thm}

\demo Clearly, distinct automorphisms of $F_n$ extend to distinct
positive automorphisms of $M$ by sending $t$ to $t$. In this sense,
we shall think $Aut(F_n)$ as a subgroup of $Aut(M)$. On the other
hand, consider $\mathbb{Z}^n \rtimes C_2=\mathbb{Z}^n \rtimes
\langle v\rangle$ so that $v^{-1}uv= -u$ for all $u\in
\mathbb{Z}^n$. It is straightforward to verify that this group acts
faithfully on $M=F_n \times \mathbb{Z}$, whereby an element
$(v^{\epsilon},u)$, $u=(u_1, \ldots ,u_n)^T$, sends $x_i$ to
$t^{u_i}x_i$ and $t$ to $t^{1-2 \epsilon}$, where $\epsilon=0,1$.
So, we shall think $\mathbb{Z}^n \rtimes C_2 \leqslant Aut(M)$.

Note that $Aut(F_n)$ and $\mathbb{Z}^n \rtimes C_2$ have trivial
intersection as subgroups of $Aut(M)$. We shall now show that they
generate $Aut(M)$. Let $\Psi \in Aut(M)$. It will be sufficient to
show that we can multiply $\Psi$ by elements in $Aut(F_n )$ and
$\mathbb{Z}^n \rtimes C_2$ until we get the identity. Note that,
since the centre of $M$ is its infinite cyclic subgroup generated by
$t$ (use $n\geqslant 2$ and see Lemma~\ref{centre}), $\langle
t\rangle$ is characteristic in $M$ and so, $t\Psi=t^{\pm 1}$. Thus,
after possibly composing with $v\in \mathbb{Z}^n \rtimes C_2$, we
can assume that $t\Psi =t$. Then, write $x_i \Psi = t^{u_i}w_i$ for
$i=1,\, \ldots ,n$, and $u=(u_1, \ldots ,u_n)^T$. Since $\Psi$ is an
automorphism, $\{w_1,\, \ldots ,\, w_n\}$ must generate (and so form
a basis of) $F_n$. Composing $\Psi$ with the automorphism which
fixes $t$ and sends $w_i$ back to $x_i$, $i=1,\ldots ,n$, we obtain
the element $(v^0, u)\in \mathbb{Z}^n \rtimes C_2$. This proves that
$Aut(F_n)$ together with $\mathbb{Z}^n \rtimes C_2$ generate
$Aut(M)$.

Consider now $\Theta \in Aut(F_n)$ and $(v^{\epsilon},u)\in
\mathbb{Z}^n \rtimes C_2$, $u=(u_1, \ldots, u_n)^T$, viewed as
elements of $Aut(M)$. We shall calculate the conjugate $\Theta
(v^{\epsilon}, u) \Theta^{-1}$ as an element of $Aut(M)$. Let
$\theta =\Theta \vert_{F_n}\in Aut(F_n)$ and let $\theta^{\rm \,
ab}=(b_{i,j})\in GL_n(\mathbb{Z})$ be its abelianisation (an
$n\times n$ integral matrix whose $i$-th row describes the total
exponent sums of $x_i \theta$). Bearing this in mind, $\Theta
(v^{\epsilon}, u) \Theta^{-1}$ acts as
 $$
\begin{array}{ccccccc}
x_i & \mapsto & x_i \Theta & \mapsto & t^{c_i}(x_i \Theta) & \mapsto
& t^{c_i} x_i
\\ t & \mapsto & t & \mapsto & t^{1-2\epsilon} & \mapsto & t^{1-2\epsilon},
\end{array}
 $$
where $c_i=\sum_{j=1}^n b_{i,j} u_j$ is the $i$-th entry of the
column vector $\theta^{\rm \, ab} u$. In other words,
 $$
\Theta (v^{\epsilon}, u) \Theta^{-1} =(v^{\epsilon}, \theta^{\rm \,
ab}u)
 $$
for every $u=(u_1, \ldots, u_n)^T \in \mathbb{Z}^n$. This
immediately shows that $\mathbb{Z}^n \rtimes C_2$ is normal in
$Aut(M)$. Hence, $Aut(F_n \times \mathbb{Z}) \cong (\mathbb{Z}^n
\rtimes C_2) \rtimes Aut(F_n)$, where the action of $Aut(F_n)$ in
this last semi-direct product is the trivial one over the $C_2$
part, and the natural one after abelianisation over the
$\mathbb{Z}^n$ part.

Lastly, to prove the final statement note that, since $\langle
t\rangle$ is central, inner automorphisms of $F_n\times \mathbb{Z}$
are just inner automorphisms by elements of $F_n$, and all of them
fix $t$. Thus, $Inn(M)=Inn(F_n)\leqslant Aut(F_n)$ and so
$Out(M)\cong (\mathbb{Z}^n \rtimes C_2)\rtimes Out(F_n)$, where the
actions are just as before but factorised by $Inn(F_n)$. $\Box$

\section{$Aut(F_2)$ and $GL_2(\mathbb{Z})$}\label{back2}

For the rest of the paper we will be considering only the case
$n=2$. Hence, we shall avoid unnecessary subscripts by using the
letters $\{a,b\}$ as free generators of $F_2 =\langle a,b\rangle$.

In this section we will briefly review some well known facts about
$Aut(F_2)$, $Out(F_2)$ and $GL_2(\mathbb{Z})$. The abelianisation
map $F_2 \to F_2^{\rm \, ab}$ induces naturally a surjective map
$Aut(F_2) \twoheadrightarrow GL_2(\mathbb{Z})$ for which, abusing
notation, we shall write $\phi \mapsto \phi^{\rm \, ab}$. More
precisely, $\phi^{\rm \, ab}$ is the $2\times 2$ integral matrix
whose first (second) row counts the total $a$- and $b$-exponent sums
of $a\phi$ (of $b\phi$). Clearly, this is well defined for any rank,
but in rank 2 it has some special properties which make this case
easier to study; the most important of these is the following well
known result (see~\cite[chapter 4, proposition 4.5]{LS}).

\begin{thm}[Nielsen]
\label{nielsen} The kernel of the abelianisation from $Aut(F_2)$ to
$GL_2(\mathbb{Z})$ consists of precisely the inner automorphisms of
$F_2$. That is, $Out(F_2) \cong GL_2(\mathbb{Z})$.
\end{thm}

This means that, for every automorphism $\phi \in Aut(F_2)$, the
$2\times 2$ integral matrix $\phi^{\rm \, ab}$ is enough to recover
the automorphism $\phi$ up to conjugation. Since the isomorphism
type of $M_{\phi}$ (and so that of $Aut(M_{\phi})$) only depends on
$\phi$ up to conjugation, $\phi^{\rm \, ab}$ contains all the
algebraic information we may want about $M_{\phi}$ and
$Aut(M_{\phi})$.

Matrices in $GL_2(\mathbb{Z})$ can be classified according to their
eigenvalues and dynamics, often leading to useful information for
$Aut(F_2)$. This is often done in the following way.

\begin{defn}
\emph{Let $A$ be a $2\times 2$ integral invertible matrix, $A\in
GL_2(\mathbb{Z})$. If $A^2$ is the identity matrix, $A^2=I_2$, we
say that $A$ is \emph{elliptic}. Otherwise, $A$ is called
\emph{hyperbolic} if $|trace(A^2)|>2$, \emph{parabolic} if
$|trace(A^2)|=2$, and \emph{elliptic} if $|trace(A^2)|<2$.}
\end{defn}

Suppose $A\in GL_2(\mathbb{Z})$ is a hyperbolic matrix. Then, $A^2
\neq I_2$ has two real eigenvalues, $\alpha, 1/\alpha$, such that
$|\alpha|>1$. Since $A$ preserves the one dimensional eigenspaces of
$A^2$, $A$ must also have two real eigenvalues, $\beta, \pm
1/\beta$, such that $|\beta|>1$. In particular, $A$ does not have
$1$ as an eigenvalue.

Suppose $A\in GL_2(\mathbb{Z})$ is a parabolic matrix. Then, $A^2
\neq I_2$ has characteristic polynomial equal to $(x \pm 1)^2$. This
implies that $A^2$ is conjugate, in $GL_2(\mathbb{Z})$, to $\pm
\left( \smallmatrix 1 & k' \\ 0 & 1
\endsmallmatrix \right)$ for some $0\neq k'\in \mathbb{Z}$ (take a rational eigenvector of eigenvalue $\pm 1$, multiply it
by a scalar to obtain an integral vector $v$ with coprime entries,
and then extend to a basis $\{u,v\}$ of $\mathbb{Z}^2$). A simple
calculation shows now that $A$ must then be conjugate to one of the
matrices
 $$
\left( \begin{array}{rr} 1 & k \\ 0 & 1 \\ \end{array} \right)
\text{ or } \left( \begin{array}{rr} -1 & k \\ 0 & -1
\\
\end{array} \right),
 $$
for some integer $k\neq 0$. These are infinite order matrices and
the first of these will turn out to be the most challenging case to
consider.

Finally, suppose $A\in GL_2(\mathbb{Z})$ is an elliptic matrix.
Then, either $A^2=I_2$ or the characteristic polynomial of $A^2$ is
equal to $x^2+1$, $x^2 +x+1$ or $x^2-x+1$. In the first case, $A$ is
either $\pm I_2$ or is conjugate, in $GL_2(\mathbb{Z})$, to $\left(
\smallmatrix 1 & k \\ 0 & -1
\endsmallmatrix \right)$ for some $k\in \mathbb{Z}$ (by similar reasoning to that above). Otherwise, $A$ has complex conjugate roots and, in
particular, it does not have 1 as an eigenvalue.

\section{The rank 2 case: proof of Theorem~\ref{main}}\label{2}

Given $\phi\in Aut(F_2)$, we shall analyse $Aut(M_{\phi})$ and prove
Theorem~\ref{main} by following the classification of matrices in
the previous section for $\phi^{\rm \, ab}$.

First, we state the following two lemmas for later use.

\begin{lem}\label{centra}
Every non-central element in $GL_2(\mathbb{Z})$ generates a finite
index subgroup of its own centraliser. The centre consists precisely
of the matrices $\pm I_2$.
\end{lem}

\demo The first statement of the lemma follows easily from the
presentation of $GL_2(\mathbb{Z})$ as an amalgamated product,
$GL_2(\mathbb{Z})\cong D_4\ast_{D_2}D_6$, where $D_n$ is the
dihedral group of order $2n$. The statement about central elements
is elementary. $\Box$

\begin{lem}\label{useful}
Consider the automorphism $\phi$ of $F_2 =\langle a,b\rangle$ given
by $a\phi=ab^k$ and $b\phi =b$, where $k\neq 0$. Then, for every
integer $r\neq 0$ and every $w\in F_2$,
\begin{itemize}
\item[i)] $Fix \, \phi= Fix \, \phi^r =\langle aba^{-1}, b\rangle$,
\item[ii)] if $w\phi^r$ is conjugate to $w$, then $w$ is conjugate to an element fixed by $\phi$.
\end{itemize}
\end{lem}

\demo Given an arbitrary reduced word $w\in F_2$, let us split it
into \emph{pieces} each of the form $b^m$, $ab^m$, $b^m a^{-1}$ or
$ab^m a^{-1}$, where $m$ is some integer. There can be many variants
of such a splitting, but we shall use the special one defined by
putting breaking points precisely before each occurrence of $a$ and
after each occurrence of $a^{-1}$. One can easily see that this
splitting is invariant under the action of $\phi$, and that the
pieces do not interact under iterates of $\phi$. Hence, if $w\phi^r
=w$ then the corresponding pieces in the splitting of $w$ must also
be fixed by $\phi^r$, which rules out the possibilities $ab^m$ and
$b^m a^{-1}$ (because $rk\neq 0$). This proves that $Fix \, \phi^r
\leqslant \langle aba^{-1}, b\rangle$. The inclusions $\langle
aba^{-1}, b\rangle \leqslant Fix \, \phi \leqslant Fix \, \phi^r$
are obvious. This proves (i).

In order to prove (ii), note that we can assume $w$ is cyclically
reduced. If $w$ is a power of $b$ there is nothing to prove. So, we
may also assume that $w$ contains $a^{\pm 1}$. Moreover, by
inverting and cyclically permuting if necessary, we may assume that
$w$ begins with $a$. So, in this particular situation, assume that
$w\phi^r$ is conjugate to $w$. The splitting of $w$ must begin with
a piece of the form $ab^m$ or $ab^m a^{-1}$, and must end with a
piece of the form $b^m$ or $ab^m$. But this splitting is stable
under iterates of $\phi$ hence, the first and last pieces in
$w\phi^r$ will be of the corresponding same types. In particular,
$w\phi^r$ is still cyclically reduced. Thus, $w\phi^r$ must be a
cyclic permutation of $w$. Then, for a suitable $s$, $w\phi^{rs}=w$
which, by (i), means that $w$ is fixed by $\phi$. This completes the
proof.\qed

{\it Proof of Theorem~\ref{main}.} Throughout the proof, let us fix
the following notation. Let $F_2 =\langle a,b\rangle$ be a free
group of rank $n=2$, let $\phi \in Aut(F_2)$ and let $M_{\phi}=F_2
\rtimes_{\phi} \mathbb{Z}$ be the mapping torus of $\phi$. Let
$\phi^{\rm \, ab} \in GL_2 (\mathbb{Z})$ be the map induced by
$\phi$ on $F_2^{\rm \, ab}\cong \mathbb{Z}^2$, i.e. the $2\times 2$
integral matrix whose rows count the total $a$- and $b$-exponent
sums of $\{ a\phi, b\phi \}$.

First of all, note that the discussions in the previous section show
that a generic matrix $\phi^{\rm \, ab}\in GL_2(\mathbb{Z})$ fits
into one of the cases distinguished in Theorem~\ref{main}. Namely,
if $A$ is hyperbolic then it satisfies (iii), if it is parabolic it
satisfies either (iii) or (v), and if it is elliptic then it fits
into either (i), (ii) or (iv). Uniqueness is a straightforward
exercise in linear algebra.

Suppose $\phi^{\rm \, ab}=I_2$. Then, $\phi =\gamma_w$ for some
$w\in F_2$. Hence, using Lemma~\ref{semi}, $M_{\phi}\cong
M_{Id}=F_2\times \mathbb{Z}$. Now, using Theorem~\ref{prod}, we have
 $$
Out(M_{\phi})\cong Out(M_{Id})\cong (\mathbb{Z}^2 \rtimes
C_2)\rtimes GL_2(\mathbb{Z}),
 $$
where the actions are the natural ones described above. This proves
Theorem~\ref{main}~(i).

Suppose $\phi^{\rm \, ab}=-I_2$. Then it does not have 1 as an
eigenvalue. Hence, by Theorem~\ref{evector}, $Out^+(M_{\phi})$ is a
normal subgroup of $Out(M_{\phi})$ of index at most two which is
isomorphic to $C(\phi^{\rm \, ab})/\langle \phi^{\rm \, ab}\rangle$,
where $C(\phi^{\rm \, ab})$ is the centraliser of $\phi^{\rm \, ab}$
in $Out(F_2 )=GL_2(\mathbb{Z})$. But $\phi^{\rm \, ab}=-I_2$, which
is central in $GL_2(\mathbb{Z})$ so,
 $$
Out^+(M_{\phi}) \cong C(\phi^{\rm \, ab})/\langle \phi^{\rm \,
ab}\rangle = GL_2(\mathbb{Z})/\{\pm I_2\} =PGL_2(\mathbb{Z}).
 $$
On the other hand, $a\mapsto a$, $b\mapsto b$, $t\mapsto t^{-1}$
determines a (well-defined) negative automorphism $\Upsilon$ of
$M_{\phi}$ and so, $Aut^+(M_{\phi})\unlhd_2 Aut(M_{\phi})$.
Furthermore, note that $\Upsilon$ has order two and commutes with
every $\Psi \in Aut^+(M_{\phi})$ (which has the form $a\mapsto w_1$,
$b\mapsto w_2$, $t\mapsto tw_3$ where $w_3$ is a \emph{palindrome},
$w_3^R =w_3$). Hence, $Aut(M_{\phi})\cong Aut^+(M_{\phi})\times
C_2$. Finally, since $Inn(M_{\phi})\leqslant Aut^+(M_{\phi})$, we
have
 $$
Out(M_{\phi})\cong Out^+(M_{\phi})\times C_2 \cong PGL_2
(\mathbb{Z})\times C_2.
 $$
This proves Theorem~\ref{main}~(ii).

Suppose that $\phi^{\rm \, ab}\neq -I_2$ does not have 1 as an
eigenvalue. Then, applying Theorem~\ref{evector} and
Lemma~\ref{centra}, we deduce that $Out(M_{\phi})$ is finite. This
proves Theorem~\ref{main}~(iii).

Now, suppose that $\phi^{\rm \, ab}$ is conjugate in
$GL_2(\mathbb{Z})$ to $\left( \smallmatrix 1 & k \\ 0 & -1
\endsmallmatrix \right)$ for some $k\in \mathbb{Z}$. Then, Proposition~\ref{oneell} shows that
$Out(M_{\phi})$ has an infinite cyclic subgroup of finite index.
This completes Theorem~\ref{main}~(iv).

Finally, suppose that $\phi^{\rm \, ab}$ is conjugate in
$GL_2(\mathbb{Z})$ to $\left( \smallmatrix 1 & k
\\ 0 & 1 \endsmallmatrix \right)$ for some $0\neq k\in \mathbb{Z}$. Then, Proposition~\ref{parab} will complete the
proof of Theorem~\ref{main}~(v) by showing that $Out(M_{\phi})$ also
has an infinite cyclic subgroup of finite index. \qed

\begin{prop}\label{oneell}
With the notation above, assume $\phi^{\rm \, ab}$ is conjugate in
$GL_2(\mathbb{Z})$ to $\left(\smallmatrix 1 & k \\ 0 & -1
\endsmallmatrix \right)$, where $k\in \mathbb{Z}$. Then
$Out(M_{\phi})$ has an infinite cyclic subgroup of finite index.
\end{prop}

\demo Using Lemma~\ref{semi}, we can assume $\phi^{\rm \, ab}
=\left( \smallmatrix 1 & k \\ 0 & -1 \endsmallmatrix \right)$.
Furthermore, composing $\phi$ by a suitable inner automorphism of
$F_2$, we can assume $a\phi =ab^k$ and $b\phi =b^{-1}$. So,
 $$
M_{\phi}=\langle a,b,t \mid t^{-1}at=ab^k,\, t^{-1}bt=b^{-1}\rangle.
 $$
It is straightforward to check that
 $$
\begin{array}{rcl} a & \mapsto & t  a \\ b & \mapsto & b \\ t & \mapsto & t.
\end{array}
 $$
defines an automorphism $\Psi$ of $M_{\phi}$ such that $[\Psi]\in
Out(M_{\phi})$ is an infinite order outer automorphism (because
inner automorphisms of $M_{\phi}$ leave $F_2$ invariant). Let us
prove now that the infinite cyclic subgroup $\langle [\Psi]\rangle$
has finite index in $Out(M_{\phi})$.

Clearly, $\phi^2=Id$ and so, $t^2$ is in the centre of $M_{\phi}$.
Also, $\langle a,b,t^2\rangle \cong F_2 \times \mathbb{Z}$ is an
index 2 subgroup of $M_{\phi}$. So, all those automorphisms of
$M_{\phi}$ which restrict to an automorphism of $\langle
a,b,t^2\rangle$ form a finite index subgroup of $Aut(M_{\phi})$
(since $M_{\phi}$ is finitely generated and so has finitely many
index 2 subgroups). Moreover, the centre of $\langle a,b,t^2\rangle$
is $\langle t^2\rangle$ so, any such automorphism sends $t^2$ to
$t^{\pm 2}$. Hence, all those automorphisms of $M_{\phi}$ which
restrict to an automorphism of $\langle a,b,t^2\rangle$ and fix
$t^2$ still form a finite index subgroup of $Aut(M_{\phi})$,
containing $\Psi^2$. Thus, we can confine our attention to the
subgroup $G\leqslant Aut(M_{\phi})$ consisting on these
automorphisms, and prove that $\langle [\Psi^2] \rangle$ has finite
index in $[G]\leqslant Out(M_{\phi})$.

Let $\Theta\in G$. Note that $[a,t]=a^{-1}t^{-1}at =b^k$ and
$[b,t]=b^{-1}t^{-1}bt =b^{-2}$ and hence $M'_{\phi}=\langle b^2
\rangle F_2^{\, \prime}$ if $k$ is even, and $M'_{\phi}=\langle b
\rangle F_2^{\, \prime}$ if $k$ is odd. Also, $M_{\phi}^{\rm \, ab}$
has torsion subgroup generated by $b$ if $k$ is even, and is
torsion-free if $k$ is odd. In any case, the preimage in $M_{\phi}$
of the (possibly trivial) torsion in $M_{\phi}^{\rm \, ab}$ is
$\langle b\rangle M'_{\phi}=\langle b\rangle F_2^{\, \prime}$. In
particular, this subgroup is characteristic in $M_{\phi}$ and so
$\Theta$ acts on $\langle a,b,t^2\rangle$ in the following way:
 $$
\begin{array}{rcl} a & \mapsto & t^r u \\ b & \mapsto & v \\ t^2 & \mapsto & t^2, \end{array}
 $$
where $u\in F_2$, $v\in \langle b\rangle F_2^{\, \prime}$, and $r$
is even. Since $\langle a, b, t^2 \rangle =\langle t^r u, v, t^2
\rangle =\langle u, v, t^2 \rangle$ and $t^2$ lies in the centre of
$M_{\phi}$, it follows that $\langle u,v \rangle=\langle
a,b\rangle$. Then, from this and the form of $v$, we deduce that
$u\in a^{\epsilon}\langle b \rangle F_2^{\, \prime}$ for some
$\epsilon =\pm 1$.

Consider now the automorphism $\Lambda=\Theta\Psi^{-\epsilon r}$
which acts like
 $$
\begin{array}{ccccl} a & \mapsto & t^r u & \mapsto & t^r (u\Psi^{-\epsilon r}) =x\in F_2 \\ b & \mapsto & v & \mapsto &
v\Psi^{-\epsilon r}=y \in F_2 \\ t^2 & \mapsto & t^2 & \mapsto & t^2
\end{array}
 $$
(note that $u$ has $a$-exponent sum equal to $\epsilon$ and so, the
$t$-exponent sum of $u\Psi^{-\epsilon r}$ is $-r$, showing that
$x\in F_2$; also, the $a$-exponent sum of $v$ is 0 and so $y\in F_2$
too). Writing $t\Lambda =t^sw$ and imposing $t^2$ to be fixed, we
deduce $t\Lambda =tz$ for some $z\in F_2$. Thus, $\Lambda$ is a
positive automorphism of $M_{\phi}$ and, by Proposition~\ref{inv},
$(\Lambda|_{F_2})^{\rm \, ab}$ lies in the centraliser of $\phi^{\rm
\, ab}$. But a straightforward matrix calculation shows that
$C(\phi^{\rm \, ab})$ is finite and so, $(x,y)$ takes only finitely
many values up to conjugacy in $F_2$. Since, by
Proposition~\ref{inv}, $z$ is uniquely determined by $(x,y)$,
$\Lambda$ also takes only finitely many values up to conjugacy,
while $\Theta$ runs over all $G$. In other words, $\langle [\Psi^2
]\rangle$ has finite index in $[G]$. \qed

\begin{prop}\label{parab}
With the notation above, assume $\phi^{\rm \, ab}$ is conjugate in
$GL_2(\mathbb{Z})$ to $\left(\smallmatrix 1 & k \\ 0 & 1
\endsmallmatrix \right)$, where $0\neq k\in \mathbb{Z}$. Then
$Out(M_{\phi})$ has an infinite cyclic subgroup of finite index.
\end{prop}

\demo Using Lemma~\ref{semi}, we can assume $\phi^{\rm \, ab}
=\left( \smallmatrix 1 & k \\ 0 & 1
\endsmallmatrix \right)$ and, furthermore, composing $\phi$ by a suitable inner automorphism of $F_2$, we can assume
$a\phi =ab^k$ and $b\phi =b$, $k\neq 0$. So, we have to understand
the automorphism group of the group
 $$
M_{\phi}=\langle a,b,t \mid t^{-1}at=ab^k,\, t^{-1}bt=b\rangle.
 $$
Before going into the analysis of $Aut(M_{\phi})$, note that both
relators have $a$- and $t$-exponent sums equal to zero. So, it makes
sense to talk about $a$- and $t$-exponent sums of elements in
$M_{\phi}$. Strictly speaking, the maps from $M_{\phi}$ to
$\mathbb{Z}$ killing $b$ and $t$ and sending $a$ to the generator of
$\mathbb{Z}$ (resp. killing $a$ and $b$ and sending $t$ to the
generator of $\mathbb{Z}$) are well defined surjective
homomorphisms. They count the total $a$- and $t$-exponent sums,
respectively (note that the notion of $b$-exponent sum in $M_{\phi}$
is only well defined modulo $k$).

It is straightforward to check that
 $$
\begin{array}{rcl} a & \mapsto & t  a \\ b & \mapsto & b \\ t & \mapsto & t \end{array}
 $$
defines an automorphism $\Psi$ of $M_{\phi}$ such that $[\Psi]\in
Out(M_{\phi})$ is an infinite order outer automorphism (because
inner automorphisms of $M_{\phi}$ leave $F_2$ invariant). Let us
prove now that the infinite cyclic subgroup $\langle [\Psi]\rangle$
has finite index in $Out(M_{\phi})$.

Consider the three automorphisms of $M_{\phi}$ defined on the
generators by
 $$
\begin{array}{ccccc}
\begin{array}{ccc} & \Omega & \\ a & \mapsto & a \\ b & \mapsto & b^{-1} \\ t & \mapsto & t^{-1} \end{array}
& \quad &
\begin{array}{ccc} & \Delta & \\ a & \mapsto & a^{-1} \\ b & \mapsto & b^{-1} \\ t & \mapsto & tb^{-k} \end{array}
& \quad &
\begin{array}{ccc} & \Xi & \\ a & \mapsto & ab \\ b & \mapsto & b \\ t & \mapsto & t \end{array}
\end{array}
 $$
(as above, checking that they are well-defined is a straightforward
exercise).

\noindent \emph{\textbf{Claim:}} for any given $\Theta \in
Aut(M_{\phi})$, there exists an integer $m$ and an element $g\in
M_{\phi}$ such that $\Theta \Psi^m \Gamma_g$ is equal to one of
$\Xi^i$, $\Xi^i \Omega$, $\Xi^i \Delta$ or $\Xi^i \Delta \Omega$,
for some $0\leqslant i\leqslant |k|-1$.

This automatically will imply that $\langle [\Psi]\rangle$ has
finite index in $Out(M_{\phi})$. In order to prove this claim note
that, since $Inn(M_{\phi})$ is a normal subgroup of $Aut(M_{\phi})$,
we may apply inner automorphisms at any point in the product $\Theta
\Psi^m$.

So, let $\Theta$ be an arbitrary automorphism of $M_{\phi}$, and
write normal forms for the images of generators, $a\Theta =t^p w_1$,
$b\Theta =t^l w_2^{\, \prime}$, $t\Theta =t^q w_3$, where
$w_1,w_2^{\, \prime},w_3\in F_2$, and $p,l,q\in \mathbb{Z}$. Write
also $w_2^{\, \prime} =w_2^r$, where $r\geqslant 1$ and $w_2$ is
either trivial or not a proper power. Applying $\Theta$ to the
equality $t^{-1}at=ab^k$ we get
 $$
w_3^{-1} t^{-q} t^p w_1 t^q w_3 =t^p w_1 (t^l w_2^{\, \prime})^k.
 $$
Comparing the $t$-exponent sums we immediately see that $kl=0$ and
hence $l=0$ and $w_2\neq 1$. Now, applying $\Theta$ to $t^{-1}bt=b$,
we also get
 $$
w_3^{-1}(w_2 \phi^q)^r w_3 =w_3^{-1}t^{-q} w_2^{\, \prime} t^q w_3 =
w_2^{\, \prime} =w_2^r.
 $$
Thus, $w_2\phi^q$ is conjugate to $w_2$ in $F_2$. By applying
Lemma~\ref{useful}~(ii), we obtain that $w_2$ is conjugate to an
element fixed by $\phi$, say $w_2=xv_2 x^{-1}$, where $x\in F_2$ and
$1\neq v_2 =v_2\phi$ is not a proper power. Now, $\Theta \Gamma_x$
is an automorphism of $M_{\phi}$ acting as
 $$
\begin{array}{ccccl} a & \mapsto & t^p w_1 & \mapsto & x^{-1}t^p w_1 x=t^p v_1 \\ b & \mapsto & w_2^{\,
\prime} & \mapsto & x^{-1}w_2^r x =v_2^r \\ t & \mapsto & t^q w_3 &
\mapsto & x^{-1}t^q w_3 x=t^q v_3,
\end{array}
 $$
where $v_1, v_3 \in F_2$. Since $b$ commutes with $t$, $v_2^r$ must
commute with $t^q v_3$. But $v_2$ commutes with $t$ since it is
fixed by $\phi$. Therefore $v_2$ commutes with $v_3$ and hence,
$v_3=v_2^s$ for some integer $s$. Finally, observe that $\{ t^p
v_1,\, v_2^r,\, t^q v_2^s \}$ must generate $M_{\phi}$. Thus, since
the $a$-exponent sum of $v_2$ is zero (by Lemma~\ref{useful}~(i)),
the $a$-exponent sum of $v_1$ must be $\pm 1$. So, without loss of
generality, we may assume that $\Theta$ acts as
 $$
\begin{array}{rcl} a & \mapsto & t^p v_1 \\ b & \mapsto & v_2^r \\ t & \mapsto & t^q v_2^s, \end{array}
 $$
where $v_1, v_2 \in F_2$, $p,q,r,s \in \mathbb{Z}$, $v_2$ is fixed
by $\phi$ and has $a$-exponent sum equal to zero, and $v_1$ has
$a$-exponent sum equal to $\epsilon =\pm 1$.

Now let $m=-\epsilon p$. It is straightforward to verify that
$\Theta \Psi^m$ acts in the following form,
 $$
\begin{array}{ccccl} a & \mapsto & t^p v_1 & \mapsto & t^p (v_1 \Psi^m ) =u_1\in F_2 \\ b & \mapsto & v_2^r &
\mapsto & u_2 \in F_2 \\ t & \mapsto & t^q v_2^s & \mapsto & t^q
u_3, \end{array}
 $$
where $u_1, u_2, u_3 \in F_2$. By Proposition~\ref{inv}, it follows
now that $\Theta \Psi^m$ restricts to an automorphism of $F_2$ with
signum $q=\pm 1$.

Consider now the automorphisms $\Theta \Psi^m$, $\Theta \Psi^m
\Omega$, $\Theta \Psi^m \Delta$ and $\Theta \Psi^m \Omega\Delta$.
Each of these leaves $F_2$ invariant and have signum $q$, $-q$, $q$
and $-q$, respectively. Also, the traces of the abelianisations of
their restrictions to $F_2$ are $d$, $e$, $-d$ and $-e$,
respectively, for some $d,e\in \mathbb{Z}$. So, one of these four
automorphisms, say $\Upsilon$, is positive and its restriction to
$F_2$ abelianises to a matrix with non-negative trace.

We shall show that, up to an inner automorphism, $\Upsilon$
coincides with $\Xi^i$ for some $0 \leqslant i \leqslant |k|-1$.
This will prove the claim since both $\Omega$ and $\Delta$ above
have order two in $Aut(M_{\phi})$.

Since $\Upsilon$ is a positive automorphism, Proposition~\ref{inv}
ensures us that the matrices $(\Upsilon \vert_{F_2})^{\rm \, ab}$
and $\phi^{\rm \, ab}=\left( \smallmatrix 1 & k \\ 0 & 1
\endsmallmatrix \right)$ do commute. But the centraliser of
$\phi^{\rm \, ab}$ in $GL_2(\mathbb{Z})$ is the set of matrices of
the form $\left( \smallmatrix 1 & * \\ 0 & 1 \endsmallmatrix
\right)$ or $\left( \smallmatrix -1 & * \\ 0 & -1 \endsmallmatrix
\right)$. So, since $(\Upsilon \vert_{F_2})^{\rm \, ab}$ has
non-negative trace, we deduce that, for some $z\in F_2$, $(\Upsilon
\Gamma_z) \vert_{F_2}$ acts as $a\mapsto ab^j$, $b\mapsto b$, for
some integer $j$. Write $j=i+\lambda k$ with $0\leqslant i\leqslant
|k|-1$, and notice that the inner automorphism $\Gamma_t$ of
$M_{\phi}$ acts as
 $$
\begin{array}{rcl} a & \mapsto & ab^k \\ b & \mapsto & b \\ t & \mapsto & t. \end{array}
 $$
Hence $\Upsilon \Gamma_z \Gamma_t^{-\lambda}$ agrees with $\Xi^i$ on
$F_2$. And they both are positive so, by Proposition~\ref{inv}~(ii),
they must also agree on $t$. Hence, up to a conjugation, $\Upsilon$
coincides with $\Xi^i$. This completes the proof. \qed

\section*{Acknowledgments}

The first named author is partially supported by the grant of the President of Russian Federation for young Doctors MD-326.2003.01, by the INTAS
grant N 03-51-3663, and by the Centre de Recerca Matem\`atica (CRM) at Barcelona. The second named author gratefully acknowledges the support of the
CRM and the UPC. The third named author is partially supported by DGI (Spain) through grant BFM2003-06613, and by the Generalitat de Catalunya
through grant ACI-013. The three authors thank the CRM for their hospitality during the academic course 2004-2005, while this paper has been written.

\end{document}